\documentclass[letterpaper, 10pt, conference]{ieeeconf}

\IEEEoverridecommandlockouts 

\usepackage[letterpaper, top=60pt, bottom=43pt, left=48pt, right=48pt]{geometry}

\usepackage{amsmath, amsfonts}
\usepackage{cancel}
\usepackage{mathtools}
\usepackage[colorlinks=true, linkcolor=black]{hyperref}
\usepackage{algorithm}
\usepackage{algpseudocode}
\usepackage{booktabs}

\usepackage{graphicx}
\usepackage[normalem]{ulem}
\usepackage{xcolor}

\newcommand{\unitary}{\mathcal U}

\newcommand{\real}[0]{\mathbb{R}}
\newcommand{\vect}[1]{\boldsymbol{\mathbf #1}}

\newtheorem{proposition}{Proposition}

\usepackage{graphics} 
\usepackage{times}
\usepackage{amssymb}  

\title{ {\LARGE \bf
A Spectral Approach to Optimal Control of the \\ Fokker--Planck Equation*} \\
{\small \it Published in IEEE Control Systems Letters (L-CSS) at \url{https://ieeexplore.ieee.org/document/11015582} }
}

\author{Dante Kalise$^{1}$, Lucas M. Moschen$^{1}$, Grigorios A. Pavliotis$^{1}$ and Urbain Vaes$^{2}$
\thanks{*DK is partially supported
by the EPSRC Standard Grant EP/T024429/1. LM is supported by the CNRS-ICL PhD Studentship program. GP is partially supported by an ERC-EPSRC Frontier Research Guarantee through Grant No. EP/X038645, ERC Advanced Grant No. 247031 and a Leverhulme Trust Senior Research Fellowship, SRF$\backslash$R1$\backslash$241055.
UV is partially supported by the European Research Council (ERC) under the EU Horizon 2020 program (grant agreement No 810367) and by the Agence Nationale de la Recherche under grants ANR-21-CE40-0006 (SINEQ) and ANR-23-CE40-0027 (IPSO).
}
\thanks{$^{1}$DK, LM, and GP are with the Dept. of Mathematics, Imperial College, London, UK
        {\tt\small \{d.kalise-balza,lmm122,pav\}@ic.ac.uk}}%
\thanks{$^{2}$UV is with Inria, Paris, France        {\tt\small urbain.vaes@inria.fr}}%
\thanks{Published in IEEE Control Systems Letters (L-CSS): \url{https://ieeexplore.ieee.org/document/11015582}}%
\thanks{© 2025 IEEE. Personal use of this material is permitted. 
Permission from IEEE must be obtained for all other uses, in any current or future media, including reprinting/republishing this material for advertising or promotional purposes, creating new collective works, for resale or redistribution to servers or lists, or reuse of any copyrighted component of this work in other works.}
\thanks{For the purpose of open access, the author has applied a Creative Commons Attribution (CC BY) license to any Author Accepted Manuscript version arising.}
}
\date{\today}

\begin{document}

\maketitle
\thispagestyle{empty}
\pagestyle{empty}


\begin{abstract}

In this paper, we present a spectral optimal control framework for Fokker--Planck equations based on the standard ground state transformation that maps the Fokker--Planck operator to a Schr\"odinger operator. 
Our primary objective is to accelerate convergence toward the (unique) steady state.
To fulfill this objective, a gradient-based iterative algorithm with Pontryagin's maximum principle and the Barzilai-Borwein update is developed to compute time-dependent controls. 
Numerical experiments on two-dimensional ill-conditioned normal distributions and double-well potentials demonstrate that our approach effectively targets slow-decaying modes, thus increasing the spectral gap.

\end{abstract}

\section{Introduction}

Optimal control of diffusion processes, whose dynamics are modeled by Fokker--Planck equations (FPEs), arises in a wide range of applications, from molecular dynamics~\cite{allen2017computer} and optics~\cite{jones2015optical} to social phenomena~\cite{helbing2006analytical, castellano2009statistical, albi2017mean, pareschi2006self}.
By designing controls that act on the FPE, one can steer the probability distribution toward a desired state, accelerate convergence to equilibrium, or stabilize otherwise slow-mixing dynamics.
In this work, we focus on the linear FPE as a fundamental case study, noting that the same framework can in principle be extended to nonlinear or mean-field PDEs via linearization.

In controlling the FPE, first- and second-order optimality conditions have been derived for an optimal control problem~\cite{aronna2021first}. 
Meanwhile, infinite-horizon control strategies for overdamped Langevin dynamics in bounded domains have been proposed to stabilize the FPE and enhance convergence to the stationary measure (Gibbs)~\cite{breiten2018control}, \cite{breiten2017reduction}, both using finite-difference schemes.
In particular, the ground-state change of variables we present in \autoref{sec:fokkerplanck} is the same as~\cite{breiten2017reduction}.
This is especially relevant to sampling, where optimal Langevin samplers improve mixing by modifying the drift term of the associated SDE~\cite{lelievre2013optimal}.
A discussion of the FP control framework, including its formulation for various FPEs and its connection to dynamic programming, can be found in~\cite{Borzi_survey}.

Unlike finite-volume or finite-difference discretizations, spectral methods rely on eigenfunction expansions to approximate solutions accurately while maintaining low computational cost. 
This approach is especially advantageous in unbounded domains, where conventional grid-based methods struggle. 
In addition, spectral formulations naturally connect to Schr\"odinger-type transformations of the Fokker--Planck operator, enabling direct use of techniques from quantum mechanics and numerical linear algebra~\cite{pavliotis2014stochastic}. 
This link was successfully exploited in multiscale stochastic differential equations, demonstrating the versatility of spectral approaches~\cite{abdulle2017spectral}.

To our knowledge, no prior work has combined, in an unbounded-domain setting, (i) the ground-state (Schrödinger) transform on $\mathbb{R}^d$, (ii) spectral discretization of the FPE, (iii) open-loop Pontryagin optimal control solved with Barzilai--Borwein (BB) solver, and (iv) systematic shape-function design to control the slowest eigenmodes.
The main contributions of our paper are:
\begin{itemize}
    \item We extend the Hermite-based approach of~\cite{mohammadi2015hermite} from quadratic to smooth confining potentials, relying on the spectrum of the associated Schrödinger operator and time-dependent controls.
    \item We derive the state-adjoint optimality system for a finite-horizon cost and propose a gradient-based iterative algorithm with a BB update to compute the optimal control.
    \item We extend the bounded-domain shape-function design of~\cite{breiten2018control} to $\real^d$ to target the slow-decaying spectral modes, and we demonstrate via numerical experiments that our method accelerates convergence toward the steady state.
\end{itemize}

The remainder of this paper is structured as follows: \autoref{sec:fokkerplanck} provides a theoretical background on the Fokker--Planck equation and the Schr\"odinger operator.
\autoref{sec:control_problem} formulates the optimal control problem. 
\autoref{sec:numerical_method} presents our numerical approach, detailing the spectral discretization. 
Finally, \autoref{sec:numerics} provides numerical experiments.

\noindent{\sl Notation}. We denote by $L^2(\mathbb{R}^d)$ the space of square-integrable functions with respect to the Lebesgue measure.
For a given probability density $\mu$ on $\mathbb{R}^d$, we say that $\phi$ belongs to the weighted space $L^2(\mathbb{R}^d, \mu)$ if $\phi \sqrt{\mu} \in L^2(\real^d)$.
For a given Hilbert space $(H, \langle\cdot, \cdot \rangle)$, we denote $\|\cdot\|_{H}$ the norm in $H$ defined through the inner product of $H$.

\section{Fokker-Planck Equation and Preliminaries}\label{sec:fokkerplanck}

Consider the diffusion process $\{X_t\}_{t\ge 0}$ in $\mathbb{R}^d$ that satisfies the stochastic differential equation
\begin{equation}
  \label{eq:sde}
  dX_t = -\nabla V\bigl(X_t\bigr)\,dt + \sqrt{2\sigma}\,dW_t,
\end{equation}
where $V \colon \mathbb{R}^d \to \mathbb{R}$ is a potential function, $\sigma > 0$ is the diffusion coefficient, and $\{W_t\}_{t \geq 0}$ denotes a standard $d$-dimensional Brownian motion.  The evolution of the probability density function (PDF) $\rho(x,t)$ of $X_t$ is governed by the forward Fokker--Planck equation
\begin{equation}
  \label{eq:fokker-planck-equation}
  \partial_t \rho = \mathcal{L}^* \rho,
\end{equation}
with the initial condition $\rho(x,0) = \rho_0(x)$, which represents the Lebesgue density function of $X_0$, satisfying $\rho_0 \ge 0$, $\int_{\mathbb{R}^d} \rho_0(x)\,dx = 1$, and $\rho_0 \in C^{\infty}(\real^d)$.
The operator $\mathcal{L}^* \colon L^2(\real^d, \rho_{\infty}^{-1}) \to L^2(\real^d, \rho_{\infty}^{-1})$ is defined as \[\mathcal{L}^*\phi \coloneqq \nabla \cdot \bigl(\phi\nabla V\bigr) \;+\; \sigma\,\Delta\phi\,.\]

We assume that $V \in C^{\infty}(\mathbb{R}^d)$ and that
\[
\lim_{|x|\to\infty} V(x) = \infty,
\, \, \text{and} \, \,
e^{-V(\cdot)/\sigma} \in L^1\bigl(\mathbb{R}^d\bigr)
\]
for all $\sigma > 0$.
Under these assumptions, the PDE is well-posed (see \cite[Thm. 2]{mohammadi2015hermite}, \cite[Thm. 4.2]{pavliotis2014stochastic}, \cite[Prop. 2.1]{breiten2018control}) and admits a unique stationary density
\[
\rho_{\infty}(x) := Z^{-1}\,e^{-\,V(x)/\sigma},
\]
where $Z \coloneqq \int_{\mathbb{R}^d} e^{-\,V(x)/\sigma}\,dx$~\cite[Prop. 4.6]{pavliotis2014stochastic}.

\vspace{0.1cm}

Next, we define
\begin{equation}
  \label{eq:W-def}
  W(x) := \frac{1}{4\sigma}\, |\nabla V(x) |^{2} - \frac{1}{2}\,\Delta V(x)
\end{equation}
and assume that $\displaystyle \lim_{|x|\to\infty} W(x) = \infty.$
Then $\rho_{\infty}$ satisfies a Poincar\'e inequality with a constant $\lambda > 0$~\cite[Cor. 1.6]{MR2386063}; 
in particular, if $\rho_{0}\in L^2\bigl(\mathbb{R}^d,\rho_{\infty}^{-1}\bigr)$,
\[
\bigl\|\rho(\cdot,t) - \rho_{\infty}\bigr\|_{L^2(\real^d, \rho_{\infty}^{-1})}
\le
e^{-\lambda\sigma t}
\,
\bigl\|\rho_0 - \rho_{\infty}\bigr\|_{L^2(\real^d, \rho_{\infty}^{-1})}.
\]
Under these hypotheses, the \emph{backward} generator $\mathcal{L}$ associated with \eqref{eq:sde} is self-adjoint in $L^2\bigl(\mathbb{R}^d,\rho_{\infty}\bigr)$ and has compact resolvent~\cite[Sec 4.4]{pavliotis2014stochastic}. 

We perform the (ground state) unitary transformation to map the Fokker-Planck operator $\mathcal{L}^*$ to a {\em Schr\"odinger operator} defined in $L^2(\mathbb{R}^d)$.  
Specifically, define the {\em unitary operator}
\[
\unitary \colon L^2(\mathbb{R}^d,\rho_{\infty}^{-1}) \to L^2(\mathbb{R}^d)
\text{ by } (\unitary \phi)(x) \coloneqq \frac{\phi(x)}{\sqrt{\rho_{\infty}(x)}}.
\]
Verifying that $\mathcal U$ is an isometric isomorphism with a well-defined inverse is straightforward.

Based on the unitary operator, we define a new operator $\mathcal{H}$ on $L^2(\mathbb{R}^d)$ by
\[
\mathcal{H} := - \unitary \,\mathcal{L}^*\, \unitary^{-1},
\]
which has the form $\mathcal{H} = -\sigma\Delta + W(x)$.
The operator $\mathcal{H}$ is self-adjoint in $L^2(\mathbb{R}^d)$ and possesses a purely discrete, unbounded spectrum with eigenvalues of finite multiplicity~\cite[Thm. XIII.67]{reed1978iv}.
Moreover, standard elliptic regularity theory guarantees that its eigenfunctions are infinitely differentiable.

Under the transformation $\psi \coloneqq \unitary \rho$, the Fokker-Planck PDE transforms to the equation
\begin{equation}\label{eq:schrodinger_equation}
      \partial_t \psi = -\mathcal{H}\psi,
\end{equation}
with initial condition $\psi(\cdot,0) = \psi_0 := \unitary \rho_0$.
Equation~\eqref{eq:schrodinger_equation} resembles a Schr\"odinger equation with imaginary time.
By transforming the Fokker--Planck operator into a Schrödinger operator, we recast the problem in the standard flat $L^2(\mathbb{R}^d)$ space, where the resulting operator is self-adjoint and its spectral properties can be analyzed using well-established results from the theory of Schr\"odinger operators.

\section{Optimal Control Problem Formulation}\label{sec:control_problem}

We consider that $X_t$ is subject to an additional conservative force:
    \begin{equation*}
        dX_t = -\nabla V\bigl(X_t\bigr)\,dt - \sum_{i=1}^{m} u_i(t) \nabla \alpha_i(X_t) \, dt + \sqrt{2\sigma}\,dW_t.
    \end{equation*}
Here, $u_i\colon [0, T] \to \mathbb{R}$ are time-dependent controls in $L^{2}([0, T])$, and $\alpha_i\colon \mathbb{R}^d \to \mathbb{R}$ are shape control functions in $C^{\infty}(\real^d)$.
In practice, $-\nabla\alpha_i$ can be generated by optical near-field traps~\cite{zaman2019fokker} or by virtual guidance in crowd and robotic applications.
The associated FPE is given by
\begin{equation}
    \label{eq:fokker-planck-control}
    \partial_t \rho = \mathcal{L}^* \rho + \sum_{i=1}^m u_i(t) \mathcal{A}_i^* \rho,
\end{equation}
where, for $i=1,\dots,m$, $\mathcal{A}_i^* \phi \coloneqq \nabla \cdot \bigl(\phi \, \nabla \alpha_i(x) \bigr)$.

The objective is to determine a pair $(\rho, \vect{u})$ that minimizes the cost functional
\begin{equation}
    \label{eq:cost_functiona_rho}
    \begin{split}
    J(\rho, \vect{u}) 
    := &\frac{1}{2}\|\rho(\cdot, T) - \rho^{\dagger}\|^2_{L^2(\real^d, \rho_{\infty}^{-1})} +
    \frac{\nu}{2} \int_{0}^{T} |\vect{u}(t) |^2 \, dt \\
    &+ \frac{\kappa}{2} \int_0^T \|\rho(\cdot, t) - \hat\rho\|^2_{L^2(\real^d, \rho_{\infty}^{-1})} \, dt
\end{split}
\end{equation}
subject to the constraint that $\rho$ satisfies equation~\eqref{eq:fokker-planck-control}, with $ \vect{u}(t) \in \mathbb{R}^{m}$, $T >0 $, $\nu > 0$, $\kappa > 0$, and $\rho^\dagger \in L^2(\real^d, \rho_{\infty}^{-1})$ is a desired final distribution.
We also include a probability distribution $\hat\rho$.
We focus on the case where $\rho^\dagger = \hat{\rho} = \rho_{\infty}$.

Applying the unitary transformation $\unitary$ and defining $\mathcal{N}_i = \unitary  \mathcal{A}_i^* \unitary^{-1}$, we transform the system into the controlled (imaginary-time) Schr\"odinger equation.
The problem then reduces to finding a pair $(\psi, \vect{u})$ that minimizes 
\begin{equation}
    \label{eq:schrodinger_optimal_control}
    \begin{split}
    J(\psi, \vect{u}) := &\frac{1}{2}\|\psi(\cdot, T) - \psi^{\dagger}\|^2_{L^2(\real^d)} + \frac{\nu}{2} \int_{0}^{T} |\vect{u}(t) |^2 \, dt \\
    &+ \frac{\kappa}{2} \int_0^T \|\psi(\cdot, t) - \hat\psi\|^2_{L^2(\real^d)} \, dt,
\end{split}
\end{equation}
where $\psi^\dagger \coloneqq  \unitary \rho^\dagger$ and $\hat\psi \coloneqq \unitary \hat\rho$, subject to the controlled equation
\begin{equation}
    \label{eq:schrodinger-control}
    \partial_t \psi = - \mathcal{H} \psi + \sum_{i=1}^m u_i(t) \, \mathcal{N}_i \psi,
\end{equation}
with each $\mathcal N_i$ modulating $\psi$ via $u_i(t)$ (see Proposition 1 for its explicit form).

\begin{proposition}
    \label{prop:operator_n}
    Define $b_i \coloneqq \nabla\alpha_i \cdot \nabla\log\sqrt{\rho_\infty}$.
    Then, for all $\phi \in L^2(\real^d)$,
    \[
    \mathcal{N}_i \phi = \nabla\cdot (\phi \nabla \alpha_i) + b_i\,\phi.
    \]
    Moreover, the formal adjoint of $\mathcal{N}_i$ in $L^2(\mathbb{R}^d)$ is
    \[
    \mathcal{N}_i^* \phi = b_i\,\phi - \nabla\alpha_i \cdot \nabla\phi.
    \]
\end{proposition}

\vspace{0.2cm}

\begin{proof}
    See the Appendix (page~\pageref{proof:operator_n}).
\end{proof}

Standard results in infinite-dimensional optimal control (see, e.g., \cite[Ch. 3]{lions1971optimal}) guarantee the existence of an optimal control given the coercivity and weak lower semicontinuity of the cost functional, together with the continuity of the control-to-state mapping.
In our setting, the quadratic structure of the cost and the compactness of the resolvent of $\mathcal{H}$ ensure that a minimizing pair $(\psi^*,\vect{u}^*)$ exists in the admissible control space $L^2(0, T;\mathbb{R}^m)$.
A similar result can be found in~\cite[Thm. 1]{peirce1988optimal}. 
In this reference, the Schr\"odinger equation from quantum mechanics (i.e. in ``real-time'') is considered. However, the proof is based on the spectral analysis of the elliptic operator and it can be applied directly to our problem.

To solve problem~\eqref{eq:schrodinger_optimal_control}--\eqref{eq:schrodinger-control}, we consider the first-order optimality conditions. 
These conditions are derived by introducing a function $\varphi \in L^2(\mathbb{R}^d)$ as the adjoint variable and forming the augmented Lagrangian
\begin{equation*}
\begin{split}
\mathcal{L}(\psi,\vect{u},\varphi) \coloneqq & J(\psi,\vect{u}) + \int_0^T \Bigl\langle \varphi(\cdot, t), \partial_t \psi(\cdot,t) \Bigr\rangle_{L^2(\mathbb{R}^d)} \, dt \\
 + \int_0^T \Bigl\langle \varphi(\cdot, &t), \mathcal{H}\psi(\cdot, t) - \sum_{i=1}^m u_i(t)\,\mathcal{N}_i\psi(\cdot, t) \Bigr\rangle_{L^2(\mathbb{R}^d)}\,dt.
\end{split}
\end{equation*}
The first-order optimality conditions are then obtained by taking variations with respect to $\psi$, $\vect{u}$, and $\varphi$, and setting the corresponding gradients to zero.

The adjoint variable $\varphi$ satisfies  
\begin{equation}
    \label{eq:adjoint_equation}
    \begin{cases}
        \partial_t \varphi = \mathcal{H} \varphi - \sum_{i=1}^m u_i(t) \mathcal{N}_i^* \varphi + \kappa(\psi - \hat\psi), \\
        \varphi(x,T) = \psi^\dagger(x) - \psi(x,T),
    \end{cases}
\end{equation}
with the optimal control given by
\begin{equation}
    \label{eq:gradient}
    \nabla_{u_i} \mathcal{L} \coloneqq \nu u_i - \langle \mathcal{N}_i^* \varphi, \psi \rangle_{L^2(\real^d)} = 0,
\end{equation}
\begin{equation}
    \label{eq:optimal_control}
    u_i(t) = \frac{1}{\nu} \langle \varphi, \mathcal{N}_i \psi \rangle_{L^2(\real^d)} = \frac{1}{\nu} \langle \mathcal{N}_i^* \varphi, \psi \rangle_{L^2(\real^d)}.
\end{equation}
In equation~\eqref{eq:adjoint_equation}, we used that $\mathcal{H}$ is a self-adjoint operator.
Therefore, an optimal pair $(\psi^*, \vect u^*)$ satisfies equations~\eqref{eq:schrodinger-control}, ~\eqref{eq:adjoint_equation}, and~\eqref{eq:optimal_control}.
%
%
\section{Numerical Method: Spectral Discretization}
\label{sec:numerical_method}

Let $\{\lambda_k, e_k\}_{k \ge 0}$ denote the eigenvalues and eigenfunctions of the operator $\mathcal{H}$.
Without loss of generality, we reorder the eigenvalues as a nondecreasing sequence so that $\lambda_0 =0$ and $\lambda_1 > 0$.
The eigenfunctions $\{e_k\}_{k \ge 0}$ form a $L^2(\real^d)$-orthonormal basis with $e_0 = \sqrt{\rho_{\infty}}$ since $\rho_{\infty}$ is the steady state of equation~\eqref{eq:fokker-planck-equation}.

The solution $\psi(x,t)$ of equation \eqref{eq:schrodinger-control} is expanded as
\[
\psi(x,t) = \sum_{k=0}^{\infty} a_k(t) \, e_k(x).
\]
Substituting this expansion into the equation \eqref{eq:schrodinger-control} and projecting onto each $e_i(x)$ yields the infinite system of ordinary differential equations (ODEs)
\begin{equation}
    \label{eq:schrodinger_solution_coeffs}
    \dot{a}_i(t) = -\lambda_i \, a_i(t) + \sum_{j=1}^m u_j(t) \sum_{k=0}^{\infty} a_k(t) \langle \mathcal{N}_j e_k, e_i \rangle_{L^2(\real^d)}
\end{equation}
Since $\lambda_k \to \infty$ as $k \to \infty$, higher modes decay rapidly, which motivates the truncation of the expansion to a finite-dimensional system.

\begin{remark}
    By construction $a_0(t)=\langle\psi(\cdot,t),e_0\rangle_{L^2(\real^d)} \equiv 1$, implying that $\dot a_0(t)=0$, which agrees with
    $\mathcal N_i^{*}[\sqrt{\rho_\infty}] = 0$.
\end{remark}

Similarly, expanding $\varphi(x,t) = \sum_{k=0}^{\infty} p_k(t) \, e_k(x)$ leads to
\begin{equation}
    \label{eq:adjoint_equation_coeffs}
    \begin{split}
        \dot{p}_i(t) &= \lambda_i \, p_i(t) - \sum_{j=1}^m u_j(t) \sum_{k=0}^{\infty} p_k(t) \big\langle \mathcal{N}^*_j e_k, e_i \big\rangle_{L^2(\real^d)} \\
        &\quad+ \kappa (a_i - \hat{a}_i),
    \end{split}
\end{equation}
where $\hat{a}_i$ are the spectral coefficients corresponding to $\hat{\psi}$.
The optimal control is then given by
\begin{equation}
    u_i(t) = \frac{1}{\nu} \sum_{j,k=0}^\infty a_j(t)\, p_k(t)\, \big\langle \mathcal{N}_i^* e_k, e_j\big\rangle_{L^2(\real^d)}.
\end{equation}

Considering a truncation to $N$ eigenmodes, we define the matrices $B^i, A^i \in \mathbb{R}^{N \times N}$ by
\[
B^i_{jk} \coloneqq \langle b_i\, e_j,\, e_k \rangle_{L^2(\real^d)}, \quad A^i_{jk} \coloneqq \langle \nabla \alpha_i \cdot \nabla e_j, e_k \rangle_{L^2(\real^d)}.
\]
The resulting coupled finite-dimensional system is
\begin{equation}
    \label{eq:coupled_system}
    \begin{cases}
        \dot{\vect a}(t) &= -\Lambda \vect a(t) + \sum_{j=1}^m u_j(t) (B^j - A^j)\vect{a}(t), \\[1mm]
        \dot{\vect p}(t) &= \Lambda \vect p(t) - \sum_{j=1}^m u_j(t) \bigl(B^j - (A^j)^T\bigr)\vect{p}(t) \\[0.2px] 
        &\quad + \kappa(\vect{a}(t) - \hat{\vect{a}}),
    \end{cases}
\end{equation}
with initial and terminal conditions
\[
a_i(0) = \langle \psi_0, e_i \rangle_{L^2(\real^d)}, \quad p_i(T) =\langle \psi^\dagger, e_i \rangle_{L^2(\real^d)} - a_i(T),
\]
where $\Lambda$ is a diagonal matrix such that $\Lambda_{ii} = \lambda_{i-1}$ for $i \in \{1, \dots, N\}$ and $B^j$ matrices are symmetric.

\begin{remark}
    When $\psi^\dagger = \hat{\psi} = \sqrt{\rho_{\infty}}$, i.e., the target state is the steady state, convergence is achieved provided that $\sum_{i=1}^{N} a_i(t)^2 \to 0$, i.e. when all modes apart from the steady-state mode decay to zero. Since modes corresponding to small eigenvalues $\lambda_i$ decay more slowly, it is important to design the control to target and accelerate the decay of them.
\end{remark}

Under the regularity assumptions on $W$, classical spectral convergence theory ensures that the projection error decays rapidly with the number of modes.
In particular, the truncated expansion converges to $\psi$ in the $L^2$ norm at a rate that depends on the smoothness of $\psi$ (see, e.g., \cite{mohammadi2015hermite}).

\subsection{The Iterative Method}
\label{sec:algorithm}

To compute the optimal control numerically, we employ a gradient descent method applied to the coupled system \eqref{eq:coupled_system}.
The space-discretization of the gradient~\eqref{eq:gradient} is given by
\begin{equation}
    \label{eq:gradient_discretization}
    \nabla_{u_j} J := \nu u_j - \big\langle \vect{p}(t), (B^j - A^j)  \vect{a}(t) \big\rangle_{\mathbb{R}^N}.
\end{equation}
Using this gradient, we perform the following iterative procedure, outlined in Algorithm~\ref{alg:schrodinger_control}.

In our numerical implementation, the eigenfunctions $\{e_k\}$ and eigenvalues $\{\lambda_k\}$ of the operator $\mathcal{H}$ are computed using a finite element approach implemented in Wolfram Mathematica via the built-in function \texttt{NDEigensystem}~\cite{Mathematica} under Dirichlet boundary conditions on a sufficiently large domain, ensuring that the potential \(W\) is large.
We emphasize that the eigenpair computation is treated as a black box: while we use finite elements in 2D for its robustness, any accurate eigensolver, such as global spectral discretizations, low-rank tensor methods, or physics-informed neural networks~\cite{jin2022physics}, can be used without affecting our convergence results.

Once the eigenfunctions are computed, we use a trapezoidal rule to obtain the matrices $B^i$ and $A^i$ for $i \in \{1,\dots,m\}$.
In a full-grid (tensor-product) spectral discretization retaining $N$ modes per coordinate, the total number of basis functions grows as $N^d$ and each coupling matrix \(A^j,B^j\) has size \(N^d\times N^d\), so memory and computational costs scale like \(O(N^{2d})\).

In Algorithm~\ref{alg:schrodinger_control}, the forward and adjoint systems are integrated using a Runge--Kutta method \(5(4)\) with relative tolerance \(1\times10^{-7}\) and absolute tolerance \(1\times10^{-9}\).
The BB step size is computed following the formulation as in \cite{barzilai1988two}.
For strictly convex quadratic objectives, BB enjoys $R$-linear convergence, i.e. exponential decay of the error toward the global minimizer, and even $Q$-linear convergence, meaning each iteration shrinks the distance to the minimizer by a fixed factor less than one, whenever the Hessian's condition number satisfies $\kappa<2$~\cite{barzilai1988two, raydan1993barzilai}.
For nonconvex objectives with Lipschitz gradients and bounded level sets, BB still drives the gradient to zero~\cite{azmi2020analysis}.
In our spectral-PMP problem, the reduced cost $J$ is smooth, coercive, strongly convex in $u$, and has bounded sublevel sets. 
Hence, $\|\nabla J(u^k)\| \to 0$ even though the iterates need not decrease monotonically.
Other approaches, such as line searches with Wolfe conditions or the nonlinear conjugate-gradient with the Armijo condition, can be employed, but we chose the BB update for its simplicity and efficiency.

\begin{algorithm}[H]
\caption{Spectral control solver with reduced gradient and Barzilai--Borwein update}\label{alg:schrodinger_control}
\begin{algorithmic}[1]
\Require \(tol>0\), \(k_{\mathrm{max}}\), initial control function \( \boldsymbol{u}^0(t)\) (and \(\boldsymbol{u}^{-1}(t)\) for BB), final time \(T\), and a time grid \(t_{\mathrm{eval}}\).
\State \textbf{Initialize:} \(k \gets 0\).
\While{\(\|\nabla J(u^k)\| > \mbox{tol}\) and \(k < k_{\mathrm{max}}\)}
    \State \textbf{(1)} Compute \(\mathbf{a}(t)\) at \(t_{\mathrm{eval}}\) by solving \eqref{eq:schrodinger_solution_coeffs} with control \(\boldsymbol{u}^k\) and initial condition \(\mathbf{a}(0)\).
    \State \textbf{(2)} Compute \(\mathbf{p}(t)\) at \(t_{\mathrm{eval}}\) by solving \eqref{eq:adjoint_equation_coeffs} with controls \(\boldsymbol{u}^k\), state \(\mathbf{a}(t)\), and terminal condition \(\mathbf{p}(T)=\mathbf{a}^\dagger-\mathbf{a}(T)\).
    \State \textbf{(3)} Evaluate \(\nabla J_{u_j}(u^k_j)\) from \eqref{eq:gradient_discretization} for each $j$.
    \State \textbf{(4)} Compute the Barzilai--Borwein step size
    \[
    \gamma_k = \frac{\sum_{j=1}^m \Bigl\langle u_j^k-u_j^{k-1},\,\nabla J_{u_j}(u_j^k)-\nabla J_{u_j}(u_j^{k-1})\Bigr\rangle}{\sum_{j=1}^m \Bigl|\nabla J_{u_j}(u_j^k)-\nabla J_{u_j}(u_j^{k-1})\Bigr|^2}.
    \]
    \State \textbf{(5)} For each \(j\in\{1,\dots,m\}\), update the control:
    \[
    u_j^{k+1}(t)=u_j^k(t)-\gamma_k\,\nabla J_{u_j}(u_j^k(t)).
    \]
    \State \textbf{(6)} Update the control interpolants on \(t_{\mathrm{eval}}\).
    \State \textbf{(7)} Set \(k \gets k+1\).
\EndWhile
\State \Return Optimal control \(\boldsymbol{u}^{k}(t)\), state coefficients \(\mathbf{a}(t)\), and adjoint coefficients \(\mathbf{p}(t)\).
\end{algorithmic}
\end{algorithm}

\subsection{Choice of the Shape Function}

In our optimal control framework, we only optimize over the time-dependent controls $u_i(t)$ while keeping the shape functions $\alpha_i(x)$ fixed.
These shape functions modulate the spatial distribution of the control input and influence the spectral properties of the operator $\mathcal{H} - u_i\mathcal{N}_i$.
In particular, a careful choice of $\alpha_i$ can increase the spectral gap of the controlled operator, thereby accelerating the convergence of the system toward its steady state.

For the purpose of choosing the functions $\alpha_i$, we introduce the deviation variable $y(x,t) = \psi(x,t) - \sqrt{\rho_{\infty}(x)}$ for every $(x,t) \in \real^d \times \real_{\ge 0}$.
The deviation function $y(x,t)$ satisfies
\begin{equation}
    \label{eq:centralized_equation}
    y_t = -\mathcal{H}y + \sum_{j=1}^m u_j \mathcal{N}_j y + \sum_{j=1}^m u_j \mathcal{B}_j,
\end{equation}
where $\mathcal{B}_j \coloneqq \mathcal{N}_j[\sqrt{\rho_{\infty}}]$ represents the effect of the control on the steady state.

Motivated by the Lyapunov-based feedback strategy introduced in~\cite{breiten2018control}, we select the shape function $\alpha_j \in L^2(\real^d)$ to align the control with the slowest decaying modes.
Specifically, we choose $\alpha_j$ such that
\[
\mathcal{B}_j = e_j, \text{ i.e. } \nabla \cdot (\rho_{\infty} \nabla \alpha_j) = \sqrt{\rho_{\infty}} e_j,
\]
for $j \in \{1, \dots, m\}$ the $m$ slowest modes of the operator $\mathcal{H}$.
This condition intuitively ensures that the control acts in the direction of the slowest decaying modes, which are typically the bottleneck for convergence.
The above PDE has a unique solution in $L^2(\real^d)$.
Moreover, by the smoothness of $\rho_{\infty}$ and~$e_j$, the solution $\alpha_j$ is smooth as well~\cite{gilbarg1977elliptic}.

Although we do not claim that this choice is optimal in any formal sense, it provides a practical means of enhancing the influence of the control on the dynamics.
In practice, we approximate the solution of this PDE by projecting on a finite-dimensional basis, that is, $\alpha_j(x) \approx \sum_{k=1}^N c_{j,k} e_k(x)$, and then solve the resulting linear system for $c_{j,k}$.

\subsection{Initialization of the Control Functions}
\label{sec:initialization}

Due to the inherent nonconvexity of the optimal control problem -- both in its infinite-dimensional formulation and its finite-dimensional discretization -- gradient-based methods are not guaranteed to find the global minimum. Such methods may become trapped in suboptimal local minima and exhibit slow convergence.
In addition, initializing the control functions as constant zero functions (i.e., $u_i \equiv 0$ for all~$i$) may yield a terminal state that is close to the desired one for some parameter choices, but not necessarily accelerating convergence.
This behavior can limit the efficiency of the algorithm and not steer the iterative process toward the desired global minimum.
Therefore, a well-informed initial guess for $u_i$ shortens the optimization time and reduces the risk of converging to a suboptimal local minimum.

Our strategy for control initialization is based on the linearization of the system around the steady state. 
Specifically, we consider the linearization of equation~\eqref{eq:centralized_equation}
\[
y_t = -\mathcal{H}\, y + \sum_{j=1}^m u_j(t)\, \mathcal{B}_j
\]
and formulate an infinite-horizon cost functional
\[
J_{\infty}(y, \vect{u}) \coloneqq \frac{1}{2}\int_{0}^{\infty} \Bigl[ \kappa\, \|\,y(\cdot, t)\|^2_{L^2(\mathbb{R}^d)} + \nu\, |\mathbf{u}(t)|^2 \Bigr] dt,
\]
which admits an optimal solution of the form
\[
u^*_j(t) = -\frac{1}{\nu} \langle \Pi y(\cdot, t), \mathcal{B}_j^* \rangle_{L^2(\real^d)},
\]
where $\Pi$ is the solution of the associated Riccati equation
{\small \[
-\mathcal{H} \Pi - \Pi \mathcal{H} - \frac{1}{\nu}\Pi \left(\sum_{j=1}^m \mathcal{B}_j^*\mathcal{B}_j \right) \Pi + \kappa\,\mathcal{I} = 0.
\]}

Here, $\mathcal{I}$ is the identity operator on $L^2(\real^d)$.
This Riccati equation admits a unique solution, providing a well-defined optimal feedback law.

We compute the feedback control $u^*_j(t)$ using spectral discretization as in the previous sections.
Since we built matrices $\Lambda$, $B^j$, and $A^j$ for the algorithm, we can use them to solve the Ricatti equation, which is reasonably fast.
This computed feedback is then used as an initial guess for the iterative gradient-based optimization algorithm from \autoref{sec:algorithm}.
Our numerical experiments indicate that this initialization accelerates convergence, reduces numerical approximation errors, and helps avoid convergence to suboptimal solutions.

\section{Numerical Experiments and Results}\label{sec:numerics}

In this section, we illustrate the performance of the proposed spectral control solver (Algorithm~\ref{alg:schrodinger_control}) on two benchmark cases.
For simplicity, we fix $\sigma=1$ in all experiments.
The first example is the quadratic potential
\[
V_1(x,y) = \frac{1}{2}(ax^2 + by^2).
\]
The invariant measure for this potential is a Gaussian measure with an ill-conditioned (when $b \ll a$) covariance matrix. This potential allows the analytical computation of the eigenvalues and eigenfunctions of the operator $\mathcal{H}$, using the Metafune--Palara--Priola theorem~\cite{lelievre2013optimal}.
In fact, for nonnegative integers $m,n$, the eigenvalues are given by $\lambda_{m,n} = an + bm$.
Since the spectral gap is $\min(a,b)$, choosing a small value for either $a$ or $b$ results in slow convergence toward the steady state. We refer to such cases as ill-conditioned.
The second example is the separable double-well potential
\[
V_2(x,y) = (x^2 - c_x)^2 + (y^2 - c_y)^2.
\]
This double-well potential is {\em metastable} for small values of $\sigma$ ({\em low temperature}), meaning probability mass remains trapped near a local minimum for long times before switching.

In all experiments, we set the final time $T=5$, control regularization parameter $\nu=10^{-4}$, and running cost regularization parameter $\kappa=5$.
The solution is approximated using $N=50$ eigenmodes, and the initial condition is chosen as a truncated Gaussian distribution centered at $(-0.2,0.5)$.
For reference, in the quadratic ill-conditioned Gaussian test, the 50th eigenvalue of $\mathcal{H}$ is \(\lambda_{50}=7\), while in the double-well example it is \(\lambda_{50}\approx48\).
We ran $500$ iterations for the gradient method, by which time the norm of the gradient falls to \(10^{-4}-10^{-3}\).
All experiments were performed on a Mac with an Apple M3 processor (8 cores) and 16 GB of memory.
As a baseline in our plots, we construct an infinite-horizon linear-quadratic-regulator LQR, Riccati-based, feedback law via the method introduced in \autoref{sec:initialization}.

\begin{figure}[!htbp]
    \setlength{\abovecaptionskip}{-15pt}
    \setlength{\belowcaptionskip}{-15pt} 
    \centering
    \includegraphics[width=\linewidth]{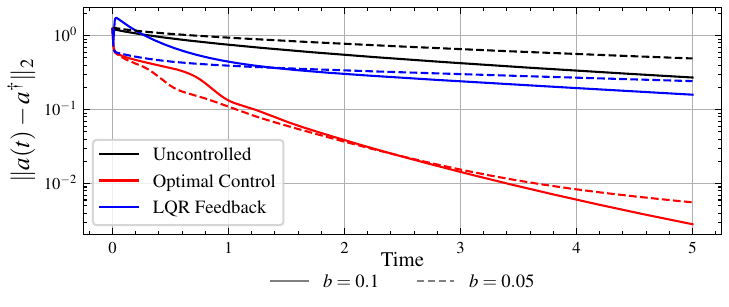}
    \caption{Time evolution of the error norm for the quadratic potential with \( a = 1 \); line styles distinguish values of \( b \) (solid for \( b = 0.1 \), dashed for \( b = 0.05 \)), while colors indicate the control strategy (black: uncontrolled, blue: LQR feedback control, red: optimal control).}
    \label{fig:convergence_v1}
\end{figure}

\begin{figure}[!htbp]
    \setlength{\abovecaptionskip}{-15pt}
    \setlength{\belowcaptionskip}{-15pt} 
    \centering
    \includegraphics[width=\linewidth]{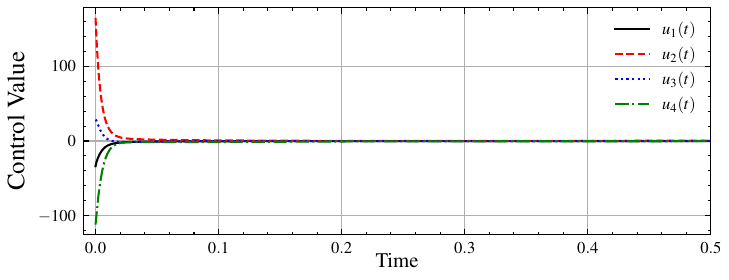}
    \caption{The four optimized control functions for the quadratic potential experiments with $a=1$ and $b=0.1$; control is shown only up to time \( t = 0.5 \) as it remains close to zero for \( t > 0.5 \).}
    \label{fig:control_profile_v1}
\end{figure}

\autoref{fig:convergence_v1} illustrates the evolution of the error norm $|a(t)-a^\dagger|_2,$ over time for the quadratic potential with $a=1$ and two choices of $b$ ($b=0.05$ and $b=0.1$).
The plot compares the optimal control (red) with the uncontrolled case (black) and the LQR feedback control (blue), demonstrating that the optimal control strategy accelerates convergence even in the challenging case of $b=0.05$.
We notice a reduction of two orders of magnitude for both experiments.
Four control functions were set to be optimized, which are depicted in \autoref{fig:control_profile_v1}.
The control signals are strongly active during the initial moments of the experiment, after which they decay as the system approaches the steady state, and the first four modes approach zero.

To assess robustness to imperfect actuator shapes, we perturbed each learned shape function by adding a Gaussian perturbation \(r_j\), i.e.\ 
\(\tilde\alpha_j(x)=\alpha_j(x)+\varepsilon\,r_j(x)\). 
We then measured the relative increase in final-time error 
\(\|\mathbf a(T)-\mathbf a^\dagger\|_2\) for 
\(\varepsilon\in\{1,5,10\}\) on the quadratic test \((a=1,\,b=0.1)\).  
The resulting percentage increases are reported in Table~\ref{tab:misalignment}.
Here, $\varepsilon$ is the $L^1$-norm of the misalignment added to each actuator.
Despite a 9.5\% increase in final-time error for a moderate misalignment ($\varepsilon = 1$), the method remains largely effective; only large perturbations cause significant performance loss.

\begin{table}[ht]
    \centering
    \begin{tabular}{c | c c c}
    \toprule
     & \(\varepsilon=1\) & \(\varepsilon=5\) & \(\varepsilon=10\) \\
    \midrule
    $\Delta(\varepsilon)$ [\%] 
      & 9.5  & 17.4   &  950 \\
    \bottomrule
    \end{tabular}
    \caption{Impact of shape-function misalignment on final-time error (quadratic test with \(b=0.1\)).\label{tab:misalignment}}
\end{table}

\begin{figure}[!htbp]
    \setlength{\abovecaptionskip}{-15pt}
    \setlength{\belowcaptionskip}{-15pt} 
    \centering
    \includegraphics[width=\linewidth]{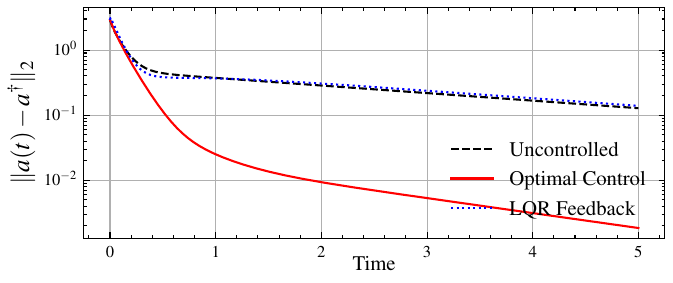}
    \caption{Time evolution of the error norm \(\|a(t)-a^\dagger\|_2\) for the double-well potential.
    The plot compares uncontrolled dynamics (dashed), LQR feedback control (dotted), and optimized open-loop control (solid).}
    \label{fig:convergence_v2}
\end{figure}

\begin{figure}[!htbp]
    \setlength{\abovecaptionskip}{-15pt}
    \setlength{\belowcaptionskip}{-15pt}
    \centering
    \includegraphics[width=\linewidth]{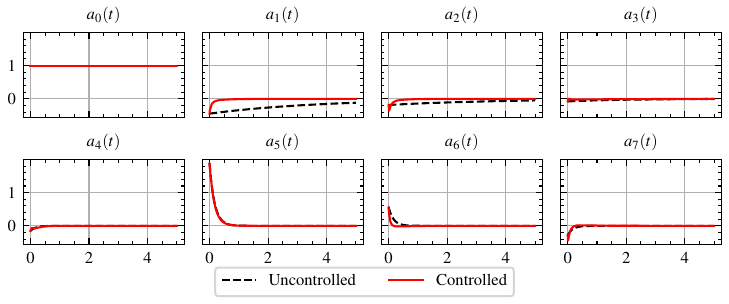}
    \caption{Time evolution of the coefficients of the spectral decomposition of the double-well potential.
    }
    \label{fig:coeffs}
\end{figure}

Figure~\ref{fig:convergence_v2} presents the time evolution of the same error for the double-well potential. 
In the uncontrolled cases, the error norm slowly decays because the system becomes temporarily trapped in one of the wells. 
In contrast, the optimal control strategy significantly accelerates convergence, as evidenced by a steeper decay of the error, finishing with a final approximation error of three orders of magnitude less.
The LQR feedback control was included as a baseline.
The two control functions optimized for this example look similar to those of \autoref{fig:control_profile_v1}.
These results confirm that our spectral control approach is robust even for multimodal potentials.
\autoref{fig:coeffs} compares the evolution of the coefficients under uncontrolled and controlled dynamics. 
In particular, the second and third components, \(a_1\) and \(a_2\), are affected by the control, whereas others (e.g., \(a_5\)) remain essentially unchanged.

\section{Conclusions}

Our numerical experiments indicate that the proposed spectral control solver effectively accelerates convergence to the desired steady state for both quadratic and double-well potentials by targeting the slower-decaying modes associated with the small eigenvalues.
Moreover, the same framework can steer an initial distribution \(\rho_0\) toward a nonequilibrium target \(\rho^\dagger\) (e.g.\ in multimodal landscapes) by choosing each $\alpha_i(x)$ to align with the dominant eigenfunctions of the linearized error operator $\psi-\psi^\dagger$.
However, because we only use time-dependent scalars $u_i(t)$ multiplying fixed spatial profiles $\alpha_i(x)$, exact controllability to $\rho^\dagger$ is not guaranteed.
Finally, as the spatial dimension increases, the required number of spectral modes grows dramatically, and the storage and manipulation of $A^j$ and $B^j$ become infeasible. 
Promising remedies include low-rank tensor decompositions, model reduction~\cite{breiten2017reduction}, sparse-grid spectral methods, or data-driven modes.  
The study of these extensions is left for future work.

\section*{APPENDIX}

We now prove \autoref{prop:operator_n}.

\begin{proof}\label{proof:operator_n}
    Applying the product rule to $\mathcal{N}_i$ yields
    \[
    \nabla\cdot(\sqrt{\rho_\infty}\, \phi\, \nabla\alpha_i)
    = \sqrt{\rho_\infty}\,\nabla\cdot(\phi\,\nabla\alpha_i)
    + \phi\,\nabla\sqrt{\rho_\infty}\cdot\nabla\alpha_i.
    \]
    Dividing by $\sqrt{\rho_\infty}$ and noting that
    \[
    \frac{\nabla\sqrt{\rho_\infty}}{\sqrt{\rho_\infty}}
    = \nabla\log\sqrt{\rho_\infty},
    \]
    we obtain $\mathcal{N}_i \phi = \nabla\cdot(\phi\,\nabla\alpha_i) + \phi\,\nabla\alpha_i\cdot\nabla\log\sqrt{\rho_\infty}$.
    This shows the desired expression.
    The formula for the formal adjoint follows from standard integration by parts.
\end{proof}

\if
The table below provides the detailed notation and main assumptions in the paper.

\begin{table}[htb]
  \centering
  \caption{Notation and main assumptions}
  \footnotesize
  \begin{tabular}{@{}ll@{}}
    \toprule
    \multicolumn{2}{l}{\it Notation}\\
    \midrule
    $\rho(x,t)$       & Fokker--Planck density \\
    $\rho_\infty(x)$  & Steady-state density \\
    $\psi=\rho/\!\sqrt{\rho_\infty}$ & Ground-state transform \\
    $\mathcal L^*$    & Forward FP operator \\
    $\mathcal H$      & Schrödinger operator, $-\sigma\Delta+W(x)$ \\
    $\mathcal N_i$    & Control operator $\mathcal U\,\mathcal A_i^*\,\mathcal U^{-1}$ \\
    $a_k(t)$          & Spectral coefficient of $\psi$ \\
    $u_i(t)$          & Time-dependent control \\
    \addlinespace[2pt]
    \multicolumn{2}{l}{\it Assumptions}\\
    \midrule
    (A1) $V\in C^\infty(\real^d),\,V(x)\to\infty$ & ensures well-posedness \\
    (A2) $e^{-V/\sigma}\in L^1(\real^d)$     & normalizable Gibbs measure \\
    (A3) $W(x)\to\infty$ as $|x|\to\infty$ & discrete spectrum of $\mathcal H$ \\
    (A4) Poincaré inequality constant $\lambda>0$ & exponential decay in $L^2(\rho_\infty^{-1})$ \\
    \bottomrule
  \end{tabular}
\end{table}
\fi

\bibliographystyle{IEEEtran}
\bibliography{biblio}

\end{document}